\documentclass{amsart}
\newtheorem{thm}{Theorem}[section]

\newtheorem{lem}[thm]{Lemma}
\newtheorem{prop}[thm]{Proposition}
\theoremstyle{example}

\theoremstyle{definition}

\theoremstyle{remark}

\numberwithin{equation}{section}

\begin{document}

\title[\bf  weighted composition operators]
{\sc\bf weighted composition operators on weak vector-valued weighted bergman spaces and Hardy spaces}

\author{Mostafa Hassanlou and Hamid Vaezi}




\begin{abstract}
In this paper we investigate weighted composition operators between weak and strong vector-valued weighted Bergman spaces and Hardy spaces.

\end{abstract}
\maketitle

\section{ \sc\bf Introduction and Preliminaries}
\indent
Weighted composition operators have been studied on different spaces of analytic functions. In \cite{cont}, Contreras and Hernandez-Diaz  have made a study of weighted composition operators on Hardy spaces whereas Mirzakarimi and Siddighi \cite{mir} have
studied these operators on Bergman and Dirichlet spaces. On Bloch-type
spaces, these operators are explored by MacCluer and Zhao \cite{mac}, Ohno \cite{ohno1},
Ohno and Zhao \cite{ohno2} and Ohno, Stroethoff and Zhao \cite{ohno3}.  In \cite{kumar} Kumar  studied  weighted composition operators between spaces of Dirichlet type by using Carleson measures. \\
 \indent
Recently these studies are about spaces of vector-valued  analytic functions.  For example, in \cite{moafa}, Wang presented some necessary and sufficient conditions for weighted composition operators to be bounded on vector-valued Dirichlet spaces and Laitila, Tylli and Wang \cite{laitila1} studied composition operators from weak to strong vector-valued Bergman spaces Hardy spaces. For some information about vector-valued Bergman spaces see \cite{ar, oscar}. \\
\indent
Let $X$ be a complex Banach space and $\mathbb{D}$ be the open unit ball of $\mathbb{C}$. We consider \textit{weight} as a strictly positive bounded
continuous function $v: \mathbb{D} \rightarrow \mathbb{\mathbb{R}}^{+}$. Let $p\geq 1$ and $v$ be a weight. The vector-valued weighted Bergman space
$A_v^p (X)$ consists of all analytic functions $f : \mathbb{D} \rightarrow X$ such that
$$ ||f||_{A_v^p (X)} = \left( \int_{\mathbb{D}} ||f(z)||_X^p v(z) dA(z) \right)^{\frac{1}{p}} < \infty. $$
where $dA$ is the normalized area measure on $\mathbb{D}$. Also, the vector-valued weighted Hardy space
$H_v^p (X)$ consists of all analytic functions $f : \mathbb{D} \rightarrow X$ for which
$$ ||f||_{H_v^p (X)}  = \sup_{0 < r < 1} \left( \int_{\mathbb{T}} ||f(r\zeta)||_X^p v(r \zeta) \ dm(\zeta) \right)^{\frac{1}{p}} < \infty, $$
where $dm(\zeta)$ is the normalized Lebesgue measure on the unit circle $\mathbb{T} = \partial \mathbb{D}$. In the case $X=\mathbb{C}$, we write
$A_v^p (X)= A_v^p$ and $H_v^p (X) = H_v^p$. Also, if $v \equiv 1$, then we have $A_v^p (X) = A^p (X)$ and $H_v^p (X)= H^p (X)$.
 The following weak versions of these spaces were considered by e.g. Blasco \cite{oscar2} and Bonet, Domanski and
Lindstrom \cite{bonet1}: the weak spaces $wA_v^p (X)$ and $w H_v^p(X)$ consist of all analytic functions $f: \mathbb{D} \rightarrow X $
for which
$$ ||f||_{wA_v^p (X)} = \sup_{||x^*|| \leq 1} ||x^* \circ f||_{A_v^p}, \ \ \ \ ||f||_{wH_v^p (X)} = \sup_{||x^*|| \leq 1} ||x^* \circ f||_{H_v^p}, $$
are finite, respectively. Here $x^* \in X^*$, the dual space of $X$. \\
\indent
Let $\varphi$ be an analytic self-map of $\mathbb{D}$; that is $\varphi(\mathbb{D}) \subset \mathbb{D}$, and $u$ a scaler-valued analytic function on
$\mathbb{D}$. We can define the weighted composition operator $uC_{\varphi}$ on the space of analytic functions as follows:
$$uC_{\varphi}(f)(z)= u(z) f (\varphi(z)).$$
When $u(z) \equiv 1$, we just have the composition operator $C_{\varphi}$, defined by $C_{\varphi}(f)= f \circ \varphi$. Also if $\varphi = I$, the identity function, then we get the multiplication operator $M_u$ defined by $M_u (f)(z) = u(z)f(z)$. It is well known that for every analytic map
$\varphi: \mathbb{D} \rightarrow \mathbb{D}$, $C_{\varphi} :A^{p}(X) \rightarrow A^{p}(X) $ and $C_{\varphi} :H^{p}(X) \rightarrow H^{p}(X) $ are bounded, and
and also on  $wA^p (X), wH^p (X)$.
For complete discussion  about composition operators we refer to \cite {cown, shap}.
We consider the infinite dimensional complex Banach space $X$, since $w A^{p}(X) = A^{p}(X)$  and  $wH^p (X) = H^p (X)$, for $\alpha > -1$ and any finite dimensional Banach space $X$. \\
But for the infinite dimensional complex Banach space $X$, $A^p (X) \not = w A^p (X)$ ($H^p (X) \not = w H^p(X)$) and
$ ||.||_{wA^p(X)} $ is not equivalent to $||.||_{A^p(X)}$ on $A^p (X)$
($ ||.||_{wH^p(X)} $ is not equivalent to $||.||_{H^p(X)}$ on $H^p (X)$), see \cite{laitila1} Proposition 3.1 (\cite{laitila2} Example 15). \\
\indent
Our aim in this paper is to compute the norm of weighted composition operators between $wA_v^p (X)$ and $A_v^p (X)$, for $p\geq 2$ and also
between $wH_v^p (X)$ and $H_v^p (X)$, for $p\geq 2$, where $v$ and $v'$ are weights.  \\
Throughout the remainder of this paper, c will denote a positive constant, the exact
value of which will vary from one appearance to the next. The notation A $\approx$ B means that
there are  positive constants $c_1$ and $c_2$ such that $c_1 A \leq B \leq c_1 A$.
\section{\sc \bf Main Results}
\begin{prop} \label{g1} \rm
Let $X$ be any complex Banach spaces, $v$ be a weight of the form $v=|\mu|$, where $\mu$ is an analytic function without any zeros on $\mathbb{D}$,
$v'$ be a weight and  $1 \leq p < \infty$. Then
\begin{equation} \label{r2}
||uC_{\varphi}: wA_v^{p}(X) \longrightarrow A_{v'}^{p}(X)|| \leq \left( \int_{\mathbb{D}} \frac{|u(z)|^p v'(z)}{(1-|\varphi(z)|^2)^{2} v(\varphi(z))} \, dA(z) \right)^{1/p}, and
\end{equation}
\end{prop}
\begin{equation} \label{r1}
||uC_{\varphi}: wH_v^p (X) \longrightarrow H_{v'}^p (X)|| \leq \sup_{0 < r<1} \left( \int_{\mathbb{T}} \frac{|u(r \zeta)|^p v'(r \zeta)}{(1-|\varphi(r \zeta)|^2) v(\varphi(r \zeta))} \, dm(\zeta) \right)^{1/p}.
\end{equation}
{\sc \bf Proof}.
\ \  By Lemma 2.1 of \cite{wolf2} we have
$$ |f(z)| \leq \frac{||f||_{A_{v}^{p}}}{(1-|z|^2)^{\frac{2}{p}} v(z)^{\frac{1}{p}}}, $$
for any $f \in A_{v}^{p}$ and $z \in \mathbb{D}$. Thus, for $f \in wA_v^p (X)$, we have
$$ ||f(z)||_X^p = \sup_{||x^*|| \leq 1} |(x^* \circ f)(z)|^p \leq \frac{1}{(1-|z|^2)^{2} v(z)} \sup_{||x^*|| \leq 1} ||x^* \circ f||_{A_{v}^{p}}^p$$
 $$ = \frac{1}{(1-|z|^2)^{2} v(z)} ||f||_{w A_{v}^{p}(X)}^p.$$
Hence
\begin{align*}
 ||uC_{\varphi} f||_{A_{v'}^{p}(X)}^p = & \int_{\mathbb{D}} |u(z)|^p ||f(\varphi(z))||_X^p v'(z) \, dA(z) \\
 \leq & ||f||_{w A_{v}^{p}(X)}^p \int_{\mathbb{D}} \frac{|u(z)|^p v'(z)}{(1-|z|^2)^{2} v(\varphi(z))} \, dA(z).
\end{align*}
The proof of the theorem is complete.
$\square$ \\
For the next results we need the following Dvoretzky$^,$s well-known theorem.
\begin{lem}\label{g2} \rm
 \cite{di} Suppose that $X$ is an infinite dimensional complex Banach space. Then for any $\epsilon > 0$ and $n \in \mathbb{N}$,
there is a linear embedding $T_n : l_n^2 \rightarrow X$ such that
\begin{equation} \label{r3}
(1+ \epsilon)^{-1} \bigg( \sum_{j=1}^{n} |a_j|^2 \bigg)^{1/2} \leq  \bigg \| \sum_{j=1}^{n} a_j T_n e_j \bigg \|_X \leq \bigg( \sum_{j=1}^{n} |a_j|^2 \bigg)^{1/2}
\end{equation}
for any scalars $a_1, a_2, \cdots, a_n$ and some orthonormal basis $\{ e_1, \cdots, e_n \}$ of $l_n^2$.
\end{lem}
Now, we prove a lower bound for the operator $uC_{\varphi}: wA_{v}^{p}(X) \longrightarrow A_{v'}^{p}(X)$, in the case $2 \leq p < \infty$.
\begin{thm} \label{g4} \rm
Let $X$ be any complex infinite-dimensional Banach space, $v$ be a weight of the form $v=|\mu|$, where $\mu$ is an analytic function without any zeros on $\mathbb{D}$,
$v'$ be a weight and  $2 \leq p < \infty$. Then
\begin{equation} \label{r4}
||uC_{\varphi}: wA_v^{p}(X) \longrightarrow A_{v'}^{p}(X)|| \approx \left( \int_{\mathbb{D}} \frac{|u(z)|^p v'(z)}{(1-|\varphi(z)|^2)^{2} v(\varphi(z))} \, dA(z) \right)^{1/p}.
\end{equation}
\end{thm}
{\sc \bf Proof}.
\ \ We only prove there exists a positive constant $c$ such that
$$ ||uC_{\varphi}: wA_v^{p}(X) \longrightarrow A_{v'}^{p}(X)|| \geq c \left( \int_{\mathbb{D}} \frac{|u(z)|^p v'(z)}{(1-|\varphi(z)|^2)^{2} v(\varphi(z))} \, dA(z) \right)^{1/p}. $$
Suppose that $x \in X$ with $||x||=1$ and define $g: \mathbb{D} \rightarrow X$ by $g(z) = \frac{1}{\mu(z)^{\frac{1}{p}}} x$. Then $g$ is an analytic function
on $\mathbb{D}$, and $||g||_{w A_v^p (X)} =1$, so that
$$||uC_{\varphi}||^p \geq ||u g\circ \varphi||_{A_{v'}^p}^p =  \int_{\mathbb{D}} \frac{|u(z)|^p v'(z)}{v(\varphi(z))} \ dA(z). $$
Hence
$$ \int_{\{ z \in \mathbb{D}: |\varphi(z)|^2 < 1/2 \}} \frac{|u(z)|^p v'(z)}{(1-|\varphi(z)|^2)^{2} v(\varphi(z))} \, dA(z) \leq
4 \int_{\mathbb{D}} \frac{|u(z)|^p v'(z)}{v(\varphi(z))} \ dA(z) \leq 4 ||uC_{\varphi}||^p. $$
So, it will be sufficient to prove that there exists a positive constant $c$ such that
$$ ||uC_{\varphi}||^p \geq c \int_{\{ z \in \mathbb{D}: |\varphi(z)|^2 \geq 1/2 \}} \frac{|u(z)|^p v'(z)}{(1-|\varphi(z)|^2)^{2} v(\varphi(z))} \, dA(z). $$
Let $\lambda_k = k^{2/p -1/2}$, for any $n \in \mathbb{N}$, we define functions $f_n$ as follows
$$ f_n (z) = \frac{1}{\mu(z)^{\frac{1}{p}}} \sum_{k=1}^{n} \lambda_k  z^{k} T_n e_k, $$
where the linear embedding $T_n$ is the same as in Lemma \ref{g2},  $||T_n||=1 $ and $||T_n ^{-1}|| \leq (1+\epsilon)$ and
$(e_1, \cdots, e_n)$ is an orthonormal basis of $\ell_2^n$. As in the proof of Theorem 3.2 \cite{laitila1}, there exists $c > 0$ such that
for $X^*$ with $||x^*|| \leq 1$, we have
\begin{align*}
\| x^* \circ f_n \|_{A_{v}^{p}} =  & \| \frac{1}{\mu(z)^{\frac{1}{p}}} \sum_{k=1}^{n} \lambda_k z^{k} x^* T_n e_k \|_{A_{v}^{p}} \\
= & ||\sum_{k=1}^{n} \lambda_k  x^* (T_n e_k) z^{k}\|_{A^{p}} \\
\leq &  c  \left( \sum_{k=1}^{n} | x^* (T_n e_k) |^2 \right)^{1/2} \leq c.
\end{align*}
It follows that $\|  f_n \|_{w A_{v}^{p}(X)} \leq c$. Thus, Fatou's lemma implies that
\begin{align*}
\| uC_{\varphi} \|^p \geq & c^{-p} \limsup_{n \rightarrow \infty} \| uC_{\varphi} f_n \|_{ A_{v'}^{p}(X)}^p \\
= & c^{-p} \limsup_{n \rightarrow \infty} \int_{\mathbb{D}} |u(z)|^p \| \frac{1}{\mu(\varphi(z))^{\frac{1}{p}}}
\sum_{k=1}^{n} \lambda_k \varphi(z)^{k}  T_n e_k \|_X^p v'(z) \, dA(z) \\
= & c^{-p} \limsup_{n \rightarrow \infty} \int_{\mathbb{D}}  \| \sum_{k=1}^{n} \lambda_k \varphi(z)^{k}  T_n e_k \|_X^p
\frac{|u(z)|^p  v'(z)}{v(\varphi(z))}  \, dA(z) \\
\geq &  \frac{c^{-p}}{(1+\epsilon)^p} \limsup_{n \rightarrow \infty} \int_{\mathbb{D}}  \left( \sum_{k=1}^{n} k^{4/p -1} |\varphi(z)|^{2k} \right)^{p/2} \frac{|u(z)|^p  v'(z)}{v(\varphi(z))} \  dA(z) \\
= &  \frac{c^{-p}}{(1+\epsilon)^p}\int_{\mathbb{D}} \left( \sum_{k=1}^{\infty} k^{4/p -1} |\varphi(z)|^{2k} \right)^{p/2}
\frac{|u(z)|^p  v'(z)}{v(\varphi(z))} \ dA(z) \\
\geq &  \frac{c_1 c^{-p}}{(1+\epsilon)^p}  \int_{\{ z \in \mathbb{D}: |\varphi(z)|^2 \geq 1/2 \}} \frac{|u(z)|^p v'(z)}
{(1-|\varphi(z)|^2)^{2} v(\varphi(z))} \ dA(z)
\end{align*}
and the last inequality is derived by Lemma 2.3 \cite{laitila1}.
As $\epsilon > 0$  was arbitrary, we obtain the desired lower bound estimate. \\
$\square$
\begin{thm} \label{g5} \rm
Let $X$ be any complex infinite-dimensional Banach space, $v$ be a weight of the form $v=|\mu|$, where $\mu$ is an analytic function without any zeros on $\mathbb{D}$,
$v'$ be a weight and  $2 \leq p < \infty$. Then
\begin{equation} \label{r5}
||uC_{\varphi}: wH_v^{p}(X) \longrightarrow H_{v'}^{p}(X)|| \approx \left( \int_{\mathbb{T}} \frac{|u(\zeta)|^p v'(\zeta)}{(1-|\varphi(\zeta)|^2) v(\varphi(\zeta))} \, dm(\zeta) \right)^{1/p}.
\end{equation}
\end{thm}
{\sc \bf Proof}. Similar to the proof of previous theorem, we only prove that there exists $c > 0$ such that 
$$||uC_{\varphi}||^p \geq c \int_{\{ \zeta \in \mathbb{T}: |\varphi(r \zeta)|^2 \geq 1/2 \}} \frac{|u(r\zeta)|^p v'(r\zeta)}{(1-|\varphi(r\zeta)|^2) v(\varphi(r\zeta))} \, dm(\zeta) .$$
Let $\lambda_k = k^{1/p - 1/2}$ and define 
$$ f_n (z) := \frac{1}{\mu(z)^{\frac{1}{p}}} \sum_{k=1}^{n} \lambda_k  z^{k} T_n e_k, $$
where the linear embedding $T_n$ is the same as in Lemma \ref{g2},  $||T_n||=1 $ and $||T_n ^{-1}|| \leq (1+\epsilon)$ and
$(e_1, \cdots, e_n)$ is an orthonormal basis of $\ell_2^n$. As in the proof of Theorem 2.2 \cite{laitila1}, there exists $c > 0$ such that
for $X^*$ with $||x^*|| \leq 1$, we have
\begin{align*}
\| x^* \circ f_n \|_{H_{v}^{p}} =  & \| \frac{1}{\mu(z)^{\frac{1}{p}}} \sum_{k=1}^{n} \lambda_k z^{k} x^* T_n e_k \|_{H_{v}^{p}} \\
= & ||\sum_{k=1}^{n} \lambda_k  x^* (T_n e_k) z^{k}\|_{H^{p}} \\
\leq &  c  \left( \sum_{k=1}^{n} | x^* (T_n e_k) |^2 \right)^{1/2} \leq c.
\end{align*}
Thus $\|  f_n \|_{w H_{v}^{p}(X)} \leq c$ and by suing  Fatou's lemma and Lemma 2.3 \cite{laitila1}, we have
\begin{align*}
\| uC_{\varphi} \|^p \geq & c^{-p} \limsup_{n \rightarrow \infty} \| uC_{\varphi} f_n \|_{ H_{v'}^{p}(X)}^p \\
= & c^{-p} \limsup_{n \rightarrow \infty} \int_{\mathbb{T}} |u(r\zeta)|^p \| \frac{1}{\mu(\varphi(r\zeta))^{\frac{1}{p}}}
\sum_{k=1}^{n} \lambda_k \varphi(r\zeta)^{k}  T_n e_k \|_X^p v'(r\zeta) \, dm(\zeta) \\
= & c^{-p} \limsup_{n \rightarrow \infty} \int_{\mathbb{T}}  \| \sum_{k=1}^{n} \lambda_k \varphi(r \zeta)^{k}  T_n e_k \|_X^p
\frac{|u(r\zeta)|^p  v'(r\zeta)}{v(\varphi(r\zeta))}  \, dm(\zeta) \\
\geq &  \frac{c^{-p}}{(1+\epsilon)^p} \limsup_{n \rightarrow \infty} \int_{\mathbb{T}}  \left( \sum_{k=1}^{n} k^{2/p -1} |\varphi(r\zeta)|^{2k} \right)^{p/2} \frac{|u(r\zeta)|^p  v'(r\zeta)}{v(\varphi(r\zeta))} \  dm(\zeta) \\
= &  \frac{c^{-p}}{(1+\epsilon)^p}\int_{\mathbb{D}} \left( \sum_{k=1}^{\infty} k^{2/p -1} |\varphi(r\zeta)|^{2k} \right)^{p/2}
\frac{|u(r\zeta)|^p  v'(r\zeta)}{v(\varphi(r\zeta))} \ dm(\zeta) \\
\geq &  \frac{c_1 c^{-p}}{(1+\epsilon)^p}  \int_{\{ z \in \mathbb{T}: |\varphi(r\zeta)|^2 \geq 1/2 \}} \frac{|u(r\zeta)|^p v'(r\zeta)}
{(1-|\varphi(r\zeta)|^2) v(\varphi(r\zeta))} \ dm(\zeta).
\end{align*}
Take $\epsilon =1$, then 
$$  \| uC_{\varphi} \|^p \geq c  \int_{\mathbb{T}} \frac{|u(r\zeta)|^p v'(r\zeta)}
{(1-|\varphi(r\zeta)|^2) v(\varphi(r\zeta))} \ dm(\zeta). $$
As $r \rightarrow 1$,
\begin{align*}
\| uC_{\varphi} \|^p \geq  & c \limsup_{r \rightarrow 1} \int_{\mathbb{T}} \frac{|u(r\zeta)|^p v'(r\zeta)}
{(1-|\varphi(r\zeta)|^2) v(\varphi(r\zeta))} \ dm(\zeta) \\
\geq  & c   \int_{\mathbb{T}} \frac{|u(\zeta)|^p v'(\zeta)}
{(1-|\varphi(\zeta)|^2) v(\varphi(\zeta))} \ dm(\zeta).
\end{align*}


\subsection*{Acknowledgment}

\end{document}